\documentclass[a4paper,11pt]{amsart}
\usepackage{amsmath, amstext, amsxtra, amsthm, amscd, amsgen,amsbsy, amsopn, amsfonts, latexsym, amssymb, amscd, amsthm}
\usepackage{marvosym}
\usepackage{color}
\usepackage{todonotes}
\usepackage{tikz-cd}
\usepackage[backend=bibtex,style=alphabetic]{biblatex}
\usepackage{nomencl}

\makenomenclature

\def\Re{\text{Re}}
\def\proofend{\hbox to 1em{\hss}\hfill $\blacksquare $\bigskip }

\newtheorem*{theorem*}{Theorem}
\newtheorem{theoremA}{Theorem}
\newtheorem{theorem}{Theorem}[section]
\newtheorem{proposition}[theorem]{Proposition}
\newtheorem{lemma}[theorem]{Lemma}
\newtheorem{remark}[theorem]{Remark}
\newtheorem*{remark*}{Remark}

\newtheorem{corollary}[theorem]{Corollary}

\newtheorem{condition}[theorem]{Condition}

\def\Z{{\mathbb Z}}
\def\R{{\mathbb R}}
\def\Q{{\mathbb Q}}
\def\C{{\mathbb C}}

\begin{document}

\title[Moduli space of sphere bundles over spheres and quotients]{Moduli space of non-negative sectional or positive Ricci curvature metrics on sphere bundles over spheres and their quotients}

\author{Jonathan Wermelinger}


\date{\today}

\maketitle
\begin{abstract}\noindent We show that the moduli space of positive Ricci curvature metrics on all the total spaces of $S^7$-bundles over $S^8$ which are rational homology spheres has infinitely many path components. Furthermore, we carry out the diffeomorphism classification of quotients of Milnor spheres by a certain involution and show that the moduli space of metrics of non-negative sectional on them has infinitely many path components. Finally, a diffeomorphism finiteness result is obtained on quotients of Shimada spheres by the same type of involution and we show that for the types that can be expressed by an infinite family of manifolds, the moduli space of positive Ricci curvature metrics has infinitely many path components. \end{abstract}

\section{Introduction and main results}

The existence problem of positively curved metrics is of fundamental importance in today's research in Riemannian geometry. For positive scalar curvature, there are by now a lot of results (see \cite[Chapter IV]{LM89}), and recently there has been substantial progress in positive Ricci curvature (see for example \cite{CW17}), but still there are only few known examples with positive sectional curvature (see \cite{Wi07}). However, a large class of examples with non-negative sectional curvature has been produced thanks to the work of Grove and Ziller \cite{GZ00}, which Goette, Kerin and Shankar \cite{GKS20} have recently extended.

Once existence has been established, one generally is lead to ask about unicity in mathematics. For Riemannian metrics, one has to choose appropriate equivalence classes, since a convex combination of any two metrics is still a metric and we therefore always have infinitely many distinct metrics in some sense. One way to try and capture different geometries, is by considering the moduli space of metrics, which is the quotient of the space of all metrics by the action of the diffeomorphism group of the manifold (via pullback metrics). In this picture, one essentially declares two metrics equivalent if they are isometric.

The study of moduli spaces is fairly recent. We only cite a few important milestones here, the interested reader is referred to \cite{TW15} for an overview of the subject. One of the most important contributions to the field has been made by Kreck and Stolz \cite{KS93}, who have defined an invariant which allows to distinguish connected components in the moduli space of positive scalar curvature metrics of certain $(4k+3)$-dimensional closed spin manifolds. This invariant has been extensively used to exhibit examples of manifolds whose moduli space of $scal>0$ or $Ric>0$ or $sec\geq 0$ metrics has infinitely many (path) connected components (see for example \cite{KS93}, \cite{Wr11}, \cite{DKT18}, \cite{Goo20a} and \cite{Goo20b} to only name a few).

Choosing a different approach, some years ago Dessai \cite{De17} proved the following result, which was independently proved by Goodman \cite{Goo20a} using the Kreck-Stolz invariant.

\begin{theoremA}\label{mainthm:s3bundles}
\cite[Theorem 4.1]{De17}\cite[Theorem A]{Goo20a} Let $M^{7}$ be the total space of a linear $S^3$-bundle over $S^4$ and assume $M^{7}$ is a rational homology sphere. The moduli space of non-negative sectional curvature metrics on $M$ has infinitely many path components. The same is true for the moduli space of positive Ricci curvature metrics.
\end{theoremA}

This result includes the so-called Milnor spheres, which are total spaces of $S^3$-bundles over $S^4$ that are homotopy 7-spheres.

Dessai uses a relative index invariant of Gromov and Lawson \cite{GL83}, which in some sense is more elementary than the Kreck-Solz invariant (which makes use of the index theory of Atiyah-Patodi-Singer). The idea of the proof is to exhibit an infinite family of diffeomorphic manifolds, by means of the diffeomorphism classification of $S^3$-bundles over $S^4$ which is due to Crowley and Escher \cite{CE03}. Using the work of Grove and Ziller \cite{GZ00}, one can construct appropriate metrics of $sec\geq 0$ on these bundles and their corresponding disk bundles. One then assumes by contradiction that the equivalence classes of such metrics with different indices are connected by a path in the moduli space of non-negative sectional curvature. It can subsequently be deduced that there is a path in the space of positive scalar curvature metrics connecting these metrics and index theory applications to positive scalar curvature finally imply a contradiction with the explicit computation of certain characteristic numbers. We will include this proof for the sake of completeness.

In an effort to produce new examples for which one can say something about the topology of the moduli space, we apply Dessai's method to $S^7$-bundles over $S^8$. It is still an open question whether these spaces admit metrics of non-negative sectional curvature, but it is straightforward to obtain positive Ricci curvature metrics on them. The diffeomorphism classification is due to Grey \cite{Gre12}.

\begin{theoremA}\label{mainthm:s7bundles}
Let $M^{15}$ be the total space of a linear $S^7$-bundle over $S^8$ and assume $M^{15}$ is a rational homology sphere. The moduli space of positive Ricci curvature metrics on $M$ has infinitely many path components.
\end{theoremA}

The next idea is to study the following spaces. Consider the involution on a Milnor sphere which is induced by fiberwise antipodal maps on $S^3$. The quotient of the Milnor sphere under this involution is homotopy equivalent to $\mathbb{R}\text{P}^7$ and will be called a Milnor projective space. The Grove-Ziller metric on the Milnor sphere is invariant under this involution and so we get a metric of $sec\geq 0$ on the quotient. We carry out the diffeomorphism classification of the Milnor projective spaces to prove the following.

\begin{theoremA}\label{mainthm:diffclassmilproj}
There are 16 different oriented diffeomorphism types of Milnor projective spaces.
\end{theoremA}

With the use of this classification, we can then get the following result about the moduli space of $sec\geq 0$ metrics on Milnor projective spaces.

\begin{theoremA}\label{mainthm:milnprojsp}
The moduli space of metrics of non-negative sectional curvature of all 16 Milnor projective spaces has infinitely many path components. The same is true for the moduli space of positive Ricci curvature metrics.
\end{theoremA}

Finally, we consider the same type of involutions on Shimada spheres (which are total spaces of $S^7$-bundles over $S^8$ that are homotopy 15-spheres). The quotients will be called Shimada projective spaces and a diffeomorphism finiteness result is obtained to prove the following. 

\begin{theoremA}\label{mainthm:shimprojsp}
There exist at least 4096 pairwise non-diffeomorphic Shimada projective spaces whose moduli space of positive Ricci curvature metrics has infinitely many path components.
\end{theoremA}

This note is organized as follows. In \textsection \ref{chapter:spherebundles}, we define sphere bundles over spheres in dimension 7 and 15 and present some of their topological properties. The Milnor and Shimada projective spaces are then defined. In \textsection \ref{sec:diffclasssphbund}, we introduce the Eells-Kuiper invariant and present the diffeomorphism classification of these sphere bundles. The diffeomorphism classification of the quotients is carried out in \textsection \ref{sec:diffclassmilshimprojspace}. The classification is based on a result by Lopez de Medrano \cite{LdM67}. We define and determine the Browder-Livesay invariant of the involution on the Milnor and Shimada spheres, as well as the normal invariants of Milnor projective spaces. Then, the Eells-Kuiper invariant of the Milnor and Shimada projective spaces are computed and used to complete the diffeomorphism classification of the former and to give a finiteness result for the latter. In \textsection \ref{sec:metricsnonnegsec} we construct metrics of $sec\geq 0$ on $S^3$-bundles over $S^4$ and the Milnor projective spaces, while in \textsection \ref{sec:posriccicurv} we construct positive Ricci curvature metrics on $S^7$-bundles over $S^8$ and the Shimada projective spaces. Finally, in \textsection \ref{sec:modspaceproof} we define the moduli space of metrics and prove our main results.

\section{Sphere bundles over spheres and their quotients by involutions}\label{chapter:spherebundles}

Let $n=1,2$ and fix a generator $\alpha\in H^{4n}(S^{4n};\Z)$. We use the same notation for the images of $\alpha$ under the isomorphisms $H^{4n}(S^{4n};\Z)\cong H_{4n}(S^{4n};\Z)\cong \pi_{4n}(S^{4n})$, where in homology $\alpha$ corresponds to the fundamental class and the second isomorphism is given by the Hurewicz map (see for example \cite[Chapter 7.4]{Sp66}). Consider $S^{4n-1}$-bundles over $S^{4n}$ with structure group $SO(4n)$. Equivalence classes of such bundles are in one to one correspondence with $\pi_{4n-1}(SO(4n))\cong\Z\oplus \Z$ (see \cite[Theorem 18.5]{Ste51}). Let $\sigma:S^{4n-1}\rightarrow SO(4n)$ be defined by $$\sigma(x)y:=xy,$$ and $\rho:S^{4n-1}\rightarrow SO(4n-1)\subset SO(4n)$ by $$\rho(x)y:=xyx^{-1},$$ where $x\in S^{4n-1}$ and $y\in \R^{4n}$ are interpreted as (unit) quaternions if $n=1$ and (unit) octonions if $n=2$, with the corresponding multiplication. Then it can be shown that $\{[\sigma],[\rho]\}$ is a free generating set of $\pi_{4n-1}(SO(4n))$. Let $M^{8n-1}_{k,l}$\label{def:totspacespherebundle} be the total space of the $S^{4n-1}$-bundle over $S^{4n}$ determined by $k[\rho]+l[\sigma]\in \pi_{n-1}(SO(4n))$ and $\pi_S$ its projection map. Hence, $M^{8n-1}_{k,l}$ can be identified with the quotient
\begin{equation}\label{eq:totspasphbuquo}
D^{4n}\times S^{4n-1}\sqcup D^{4n}\times S^{4n-1}/\sim
\end{equation}
where $(x,f(x)y)\sim (x,y)\in S^{4n-1}\times S^{4n-1}$ for the clutching function $$f:S^{4n-1}\rightarrow SO(4n):x\rightarrow (y\mapsto x^{k+l}yx^{-k}).$$ Let $D^{4n}\rightarrow W^{8n}_{k,l}\xrightarrow{\pi_D} S^{4n}$ be the associated disk bundle\label{def:diskbundle} and denote by $\xi_{k,l}$ the associated vector bundle\label{def:vectorbundle} $\R^{4n}\rightarrow E_{k,l}\xrightarrow{\pi_E}S^{4n}$.

Since $W^{8n}_{k,l}\simeq S^{4n}$ we have $H^{4n}(W^{8n}_{k,l};\Z)\cong \Z$. We orient $W^{8n}_{k,l}$ in such a way that $\text{sign}(W^{8n}_{k,l})=1$ (see p.\pageref{def:signaturemfdwithboundary} for the defintion) and fix the induced orientation on the boundary $M^{8n-1}_{k,l}$.

\begin{remark}\label{rm:diffeo}
Note that a change of orientation of the base leads to a diffeomorphism $M^{8n-1}_{k,l}\cong- M^{8n-1}_{-k,-l}$, whereas a change of orientation in the fiber leads to $M^{8n-1}_{k,l}\cong- M^{8n-1}_{-k-l,l}$. Hence $M^{8n-1}_{k,-l}\cong M^{8n-1}_{k-l,l}$ and we can therefore focus on $l\geq 0$ from now on.
\end{remark}

We summarize some properties of these bundles and spaces in the following (see \cite{CE03} and \cite{Gre12}).

\begin{theorem}\label{thm:spherebundleproperties}\hfill
\begin{enumerate}
\item The Euler class of $\xi_{k,l}$ is $e(\xi_{k,l})=l\alpha \in H^{4n}(S^{4n};\Z)$.
\item The integer cohomology groups of $M^{8n-1}_{k,l}$ are
$$H^0(M^{8n-1}_{k,l};\Z)\cong H^{8n-1}(M^{8n-1}_{k,l};\Z)\cong \Z,$$ 
$$H^{4n}(M^{8n-1}_{k,l};\Z)\cong \Z_l,$$
$$H^{j}(M^{8n-1}_{k,l};\Z)=0 \text{ otherwise}.$$
\item Both $W^{8n}_{k,l}$ and $M^{8n-1}_{k,l}$ are spin and both have a unique $Spin$ structure.
\item The only non-trivial Pontrjagin classes of $\xi_{k,l}$, $W^{8n}_{k,l}$ and $M^{8n-1}_{k,l}$ are
$$p_{n}(\xi_{k,l})= (4n-2)(2k+l) \alpha\in H^{4n}(S^{4n};\Z)$$
 $$p_{n}(W^{8n}_{k,l})= (4n-2)(2k+l) \pi^*_D(\alpha)\in H^{4n}(W^{8n}_{k,l};\Z)$$
 $$p_{n}(M^{8n-1}_{k,l})= (4n-2)2k \pi^*_S(\alpha)\in H^{4n}(M^{8n-1}_{k,l};\Z)$$ respectively.

\end{enumerate}
\end{theorem}

Since $H^{4n}(M^{8n-1}_{k,l};\Z)\cong \Z_l,$ it follows that $M^{8n-1}_{k,l}$ and $M^{8n-1}_{k',l'}$ cannot be homotopy equivalent if $l\neq l'$ (hence they cannot be diffeomorphic).

If $l=0$, $M^{8n-1}_{k,0}$ is not a rational homology sphere. This case will be excluded from now on.

If $l=1$, it was proved by Milnor \cite{Mi56} for $n=1$ and by Shimada \cite{Sh57} for $n=2$ that $M^{8n-1}_{k}:=M^{8n-1}_{k,1}$ is homeomorphic, but not always diffeomorphic, to the standard $(8n-1)$-sphere\footnote{Note that Milnor \cite{Mi56} and Shimada \cite{Sh57} use different generators of $\pi_{4n-1}(SO(4n))$.}. Consequently, $M^{8n-1}_k$ will be called a \emph{Milnor sphere} if $n=1$ and a \emph{Shimada sphere} if $n=2$.

Now consider the involution $\tau$ on $M^{8n-1}_{k}$ which is induced by the fiberwise antipodal map on $S^{4n}$. Indeed, the antipodal map commutes with the action of the structure group $SO(4n)$ on the fibers and thus induces an action on the total space $M^{8n-1}_k$ (see \cite[II.1.1]{Bre72}). Equivalently, $\tau$ is the map induced by $(x,y)\mapsto (x,-y)$ on $(x,y)\in D^{4n}\times S^{4n-1}$ from Equation (\ref{eq:totspasphbuquo}), when descending to the quotient. For $n=1$, the pair $(M^7_k,\tau)$ is called a \emph{Hirsch-Milnor involution}. For both $n=1,2$, this involution is smooth, orientation preserving and fixed-point free. The quotient space $Q^{8n-1}_{k}:=M^{8n-1}_{k}/\tau$\label{def:milnorshimadaprojspace} is homotopy equivalent to $\R \text{P}^{8n-1}$ (see \cite[(3.1) Proposition]{Bro67}) and will be called a \emph{Milnor projective space} if $n=1$ and a \emph{Shimada projective space} if $n=2$. Since being spin is a homotopy invariant, it follows that $Q^{8n-1}_{k}$ is spin for both $n=1,2$ and all $k$ (see \cite[p.86-87]{LM89}). 

We also denote the involution induced by fiberwise antipodal maps on $W^{8n}_k$ by $\tau$. The fixed point set of this involution is the zero-section $S_0\cong S^{4n}$.

\begin{remark}\label{rmk:spinstructmilshimproj}
Suppose that $W^{8n}_k$ is equipped with a $\tau$-invariant metric which is of product form near the boundary $M^{8n-1}_k$. Since there is a unique $Spin$ structure on $W^{8n}_k$, the involution $\tau$ preserves this $Spin$ structure (and the one on $M^{8n-1}_k$). The fixed point set has even codimension in $W^{8n}_k$, and so by \cite[Proposition 8.46]{AtBo68} it follows in their terminology that $\tau$ is of \emph{even} type. This means that the group action induced by $\Z_2=\{Id,\tau\}$ lifts to a $\Z_2$-action on the $Spin$ structure on $W^{8n}_k$, the complex spinor bundle and its space of sections which commutes with the $Spin^+$ Dirac operator $D^+_W$ (see Appendix A for the definition). Hence, the $Spin$ structure on $M^{8n-1}_k$ descends to a $Spin$ structure on $Q^{8n-1}_{k}$ and its $Spin$ Dirac operator $D_M$ (which is the restriction of $D^+_W$ to the boundary) commutes with the induced action of $\tau$ on sections of the spinor bundle and thus descends to a $Spin$ Dirac operator $D_Q$ on $Q^{8n-1}_{k}$.
\end{remark}

\section{Diffeomorphism classification of sphere bundles over spheres}\label{sec:diffclasssphbund}

Let $M^{4k-1}$ be a closed, oriented $(4k-1)$-dimensional manifold. Let $W^{4k}$ be a compact, spin manifold with boundary $\partial W=M$. The spin structure on $W$ restricts to a spin structure on $M$.

Suppose furthermore that the the following holds.

\begin{condition}[Condition $\mu$]\label{cond:mu}
\begin{enumerate}
\item The homomorphisms

$$j^*:H^{4i}(W,M;\Q)\rightarrow H^{4i}(W;\Q) \qquad 0<i<k$$
$$j^*:H^{2k}(W,M;\Q)\rightarrow H^{2k}(W;\Q)$$
in the exact sequence of the pair $(W,M)$ are isomorphisms.
\item The homomorphism $i^*:H^1(W;\Z_2)\rightarrow H^1(M;\Z_2)$ is surjective, where $i:M\rightarrow W$ denotes the inclusion.

\end{enumerate}
\end{condition}

\noindent Under these conditions, we can define
\begin{equation}\label{eqn:overlinepontrj}
\overline{p}_i(W):=(j^*)^{-1}(p_i(W))\in H^{4i}(W,M;\Q),
\end{equation} $0<i<k$, where $p_i(W)\in H^{4i}(W;\Q)$ are the rational Pontrjagin classes of $W$.

If $M$ and $W$ satisfy Condition \ref{cond:mu}, we can define the \emph{Eells-Kuiper invariant} of $M$:

\begin{equation}\label{eqn:eellskuiperinv}
\mu(M)\equiv \frac{1}{a_k}\Big(\langle N_k(\overline{p}),[W,M]\rangle+t_k \text{sign}(W)\Big) \text{ mod }1
\end{equation}
where
$$N_k(\overline{p}):=\hat{A}_k(\overline{p}_1(W),...,\overline{p}_{k-1}(W),0)-t_k L_k(\overline{p}_1(W),...,\overline{p}_{k-1}(W),0),$$
$a_k:=4/(3+(-1)^k)$ and $t_k:=\hat{A}_k(0,...,0,1)/L_k(0,...,0,1)$. Here $\hat{A}$ and $L$ denote the respective genera from \textsection \ref{section:spindiracoperator} and \textsection \ref{section:signatureoperator}.

\begin{proposition}\label{prop:eellskuipproperties}
\cite[\textsection 3.]{EK62} The Eells-Kuiper invariant satisfies the following properties.
\begin{enumerate}
\item If $M_1$ and $M_2$ are orientation preservingly diffeomorphic, then $\mu(M_1)=\mu(M_2)$.
\item If $-M$ denotes $M$ with opposite orientation, then $\mu(-M)=-\mu(M)$.
\item $\mu(M_1\# M_2)=\mu(M_1)+\mu(M_2)$.
\end{enumerate}
\end{proposition}

Let $M=M^{8n-1}_{k,l}$ be the total space of a $S^{4n-1}$-bundles over $S^{4n}$ for $n=1,2$ and $W=W^{8n}_{k,l}$ the total space of the corresponding disk bundle. Let $x=\pi^*_D(\alpha)\in H^{4n}(W;\Z)\cong \Z$ and let $y$ be a generator of $H^{4n}(W,M;\Z)\cong \Z$ such that $j^*(y)=lx$ where $j^*:H^{4n}(W,M;\Z)\rightarrow H^{4n}(W;\Z)$ is the homomorphism from the long exact sequence of the pair. Hence $(j^*)^{-1}(x)=\frac{1}{l}y$ and from now on, we are dealing with rational coefficients. Now we have $$\langle\overline{p}^2_n(W),[W,M]\rangle=\langle (j^*)^{-1}(p_n(W))\cup p_n(W),[W,M]\rangle,$$ where $\overline{p}_n(W)=(j^*)^{-1}(p_n(W))$ (note that $j^*$ is an isomorphism on the cohomology groups with rational coefficients in degree $4n$). Using Theorem \ref{thm:spherebundleproperties}, 
we compute
\begin{equation}\label{eqn:pontrdiskbund}
\langle\overline{p}^2_n(W^{8n}_{k,l}),[W^{8n}_{k,l},M^{8n-1}_{k,l}]\rangle=(4n-2)^2\frac{(2k+l)^2}{l}.
\end{equation} Using Equation (\ref{eqn:eellskuiperinv}), we obtain the following.

\begin{lemma}\label{lemma:eellskuipspherebundles}
The Eells-Kuiper invariant of $M^{8n-1}_{k,l}$ is given by 

$$\mu(M^{8n-1}_{k,l})\equiv \frac{1}{2^{4n-2}q_n}\frac{(2k+l)^2-l}{8l} \emph{ mod }1,$$ where $q_1=7$ and $q_2=127$.
\end{lemma}

We can now present the diffeomorphism classification of the sphere bundles.

\begin{theorem}\label{thm:diffeoclassspherebun}
\cite[Theorem 1.5]{CE03}\cite[Theorem 3.8.3]{Gre12} Let $M^{8n-1}_{k,l}$ and $M^{8n-1}_{k',l}$ be the total spaces of two $S^{4n-1}$-bundles over $S^{4n}$ for $n=1,2$, $l>0$. Then $M^{8n-1}_{k,l}$ is orientation preservingly diffeomorphic to $M^{8n-1}_{k',l}$ if and only if $$\mu(M^{8n-1}_{k,l})=\mu(M^{8n-1}_{k',l}) \quad \text{and}$$ $$2k\equiv 2\gamma k'\emph{ mod } l$$ for some $\gamma$ satisfying $\gamma^2\equiv 1 \emph{ mod }l$.  
\end{theorem}

In both cases, the proof makes use of the classification of highly connected manifolds in dimensions 7 and 15 from Crowley's PhD thesis \cite{Cr02}. See \cite[Chapter 1]{Cr02} for an overview of this classification.

From the above classification theorem, we can immediately conclude the following.

\begin{corollary}\label{cor:infinitefamilysphbund}
Let $n=1,2$. For each total space $M^{8n-1}_{k,l}$ of an $S^{4n-1}$-bundle over $S^{4n}$, the set $\{M^{8n-1}_{k',l}\}_{m\in\Z}$, $k'=k+2^{4n-2}lm\cdot q_n$, $q_1=7$ and $q_2=127$, is an infinite family of manifolds all orientation preservingly diffeomorphic to $M^{8n-1}_{k,l}$.
\end{corollary}

With the help of some modular arithmetic, one can deduce the number of different diffeomorphism types in case $l=1$.

\begin{corollary}\label{cor:numberofdifftypesspherebund}
\cite[\textsection 6. and 9.]{EK62} There are 16 different oriented diffeomorphism types of Milnor spheres and 4096 different oriented diffeomorphism types of Shimada spheres.
\end{corollary}

\section{Diffeomorphism classification of quotients}\label{sec:diffclassmilshimprojspace}

The diffeomorphism classification of the quotients will be based on the following result.

\begin{theorem}\label{thm:projspbrownorm}
\cite[Theorem 4]{LdM67} Let $Q^n$ be a smooth manifold and $h:Q^n\rightarrow \R \emph{P}^n$ a homotopy equivalence, $n\geq 5$. Then the diffeomorphism class of $Q^n$ is determined, up to connected sum with an element of $bP_{n+1}$, by its Browder-Livesay invariant and its normal invariant.
\end{theorem}

Recall that if $\Theta_n$\label{def:diffgrouphomotopyspheres} is the group of $h$-cobordism classes of $n$-dimensional homotopy spheres, where group addition is by connected sum, then $bP_{n+1}\subset \Theta_n$ is the subgroup consisting of those elements which are the boundary of a parallelizable manifold.

We first define normal invariants. Let $f:X \rightarrow M$ be a map between two smooth manifolds. Let $\nu_X$ denote the stable normal bundle over $X$, $\xi$ a stable vector bundle over $M$ and $b:\nu_X\rightarrow \xi$ a bundle map covering $f$. If $b':\nu_X\rightarrow \xi'$ is another bundle map with $\xi'$ a stable vector bundle over $M$, then $b$ and $b'$ are \emph{equivalent} if there exists a bundle isomorphism $c:\xi\rightarrow\xi'$ such that $c\circ b = b'$.

With this notation, a \emph{normal map} is a pair $(f,[b])$ where $f:X^n\rightarrow M$ is a map of degree one and $[b]$ an equivalence class of bundle maps. Two normal maps $(f_i,[b_i]):X_i\rightarrow M$, $i=0,1$, are called \emph{normally cobordant} if there exists a map $F:Y^{n+1}\rightarrow M$ and a bundle map $B:\nu_Y\rightarrow \xi$ covering $F$, with $Y^{n+1}$ a cobordism between $X_0$ and $X_1$, and such that $F|_{X_i}=f_i$ as well as $[B|\nu_{X_i}]=[b_i]$ for $i=0,1$. Then the pair $(F,[B])$ is called a \emph{normal cobordism} between $(f_0,[b_0])$ and $(f_1,[b_1])$. The normal cobordism class of a manifold $M$ is called its \emph{normal invariant} and the set of normal invariants of $M$ will be denoted by $\mathcal{N}(M)$\label{def:setofnormalinvariants}.

Let $G_n:=\{f:S^{n-1}\rightarrow S^{n-1}|deg(f)=\pm 1\}$, which is a topological monoid when equipped with the compact open topology, define the direct limit $G:=\text{lim}_{n\rightarrow \infty}G_n$ via suspension and consider the corresponding classifying space $BG$. Denote by $BO$ the classifying space for stable linear bundles. Then there is a fibre map $BO\rightarrow BG$ whose fibre is denoted by $G/O$. For a smooth manifold $M$, $\mathcal{N}(M)$ is non-empty and it is in one-to-one correspondence with $[M,G/O]$\label{def:[X/GO]} (see \cite[Theorem 2.23]{MM79}). From now on, we will identify the set of normal invariants with $[M,G/O]$ without further mention.

In order to study involutions on homotopy spheres, we introduce the Browder-Livesay invariant.

Let $M^n$ be a closed oriented smooth $n$-dimensional manifold and $T:M\rightarrow M$ a smooth fixed point free involution. A \emph{characteristic submanifold} of $(M,T)$ is a compact submanifold $C^{n-1}\subset M^{n}$ such that there exists a manifold with boundary $A^n$ satisfying $C=A\cap T(A)$, $M=A\cup T(A)$ and $\partial A=C$. We will also say that $P:=C/T$ is a characteristic submanifold for the quotient $M/T$.

Let $(\Sigma^{4k+3},T)$, $k\geq 1$, be a homotopy $(4k+3)$-sphere with a smooth fixed point free involution $T$. Let $C\subset \Sigma^{4k+3}$ be a characteristic submanifold of this involution. For $x,y\in \text{ker}(H_{2k+1}(C;\Z)\rightarrow H_{2k+1}(A;\Z))/\text{torsion}$, where the map is induced by the inclusion, we consider the bilinear form $$B(x,y):=x\cdot T_*(y)$$ where the dot stands for the intersection number. The bilinear form $B$ is even, symmetric and unimodular \cite[I.1.3]{LdM71}. It follows that the index of $B$, defined as the difference between the number of positive and negative values on the diagonal of a diagonalization of $B$, is a multiple of 8 (see for example \cite[p.92 Korollar]{HM68}).

We can now define\footnote{We only give the definition in dimensions $4k+3$, but the Browder-Livesay invariant can be defined in all dimensions (see \cite[I.1.3]{LdM71}).} the \emph{Browder-Livesay invariant} by \label{def:browderlivesayinv} $$\sigma(\Sigma^{4k+3},T):=\frac{1}{8}\text{index}(B),$$ which by the above observation is an integer.
It can be shown that this invariant is well-defined (i.e. it does not depend on the choice of the characteristic submanifold, see \cite[Lemma. 3.2.]{BL73}).

We say that an involution $(\Sigma^n,T)$ \emph{desuspends} if there is a smoothly embedded $S^{n-1} \subset \Sigma^n$ such that $T(S^{n-1})=S^{n-1}$. The involution \emph{doubly desuspends} if there is also a smoothly embedded $S^{n-2}\subset S^{n-1} \subset \Sigma^n$ such that $T(S^{n-2})=S^{n-2}$. The Browder-Livesay invariant gives a condition for an involution on a homotopy sphere to desuspend.

\begin{theorem}\label{thm:desuspensionthm}
\cite[I.1.3 Theorem]{LdM71} For $n\geq 6$, a smooth fixed point free involution $(\Sigma^{n},T)$ desuspends if and only if $\sigma(\Sigma^{n},T)=0$.
\end{theorem}

\subsection{Browder-Livesay invariant of involution on Milnor and Shimada spheres}

Let $M^{8n-1}_k$ be either a Milnor or Shimada sphere and consider the involution $\tau$ induced by fiberwise antipodal maps.

We begin by computing the Browder-Livesay invariant of $(M^{8n-1}_k,\tau)$.

\begin{theorem}\label{thm:blhirschmilnor}
For each $k\in \Z$, the involution $\tau$ on $M^{8n-1}_k$ doubly desuspends if $n=1$ and desuspends if $n=2$. Hence, in particular $\sigma(M^{8n-1}_k,\tau)=0$ for both $n=1,2$.
\begin{proof}

For $n=1$, the proof is due to Hirsch and Milnor (see \cite[Lemma 1]{HM64}). The same argument applies to $n=2$.

%
%
%

We use the following explicit description of $M^{15}_k$ by Shimada \cite{Sh57}. Let $S^8=\{(s,\sigma)\in \mathbb{O}\times \R \,|\, \|s\|^2+(\sigma -\frac{1}{2})^2=\frac{1}{4}, 0\leq \sigma \leq 1\}\subset\R^9$, where $\mathbb{O}\cong \R^8$ denotes the octonions. Let $V_1=S^8\setminus{\{(0,0)\}}$ and $V_0=S^8\setminus{\{(0,1)\}}$. Then
$$M^{15}_k=V_1\times S^7\cup_\psi V_0\times S^7$$ where $\psi:(V_1\cap V_0)\times S^7\rightarrow (V_1\cap V_0)\times S^7$ is the diffeomorphism
$$\psi((s,\sigma,t)_1)=\Big(s,\sigma,\frac{s^{k+1}ts^{-k}}{\|s\|}\Big)=(s,\sigma,t')_0.$$ Define $h:M^{15}_k\rightarrow \R$ by
$$h([s,\sigma,t])=\sqrt{\sigma}\Re(t), \qquad h([s,\sigma,t'])=\frac{\Re(\overline{s}t')}{\sqrt{1-\sigma}}.$$ Then $h$ has two non-degenerate critical points $(0,1,\pm 1)$. Therefore, it follows from Morse theory that $S^{14}_0:=h^{-1}(0)$ is diffeomorphic to the standard 14-sphere. It is easy to see that $S^{14}_0$ is invariant under $\tau$ and therefore this involution desuspends.

The last statement follows by Theorem \ref{thm:desuspensionthm}.

\end{proof}
\end{theorem}

\subsection{Normal invariants of Milnor projective spaces}

Next, we determine the normal invarinats of Milnor projective spaces $Q^{7}_k:=M^{7}_k/\tau$.

Let $M$ be a smooth manifold. The \emph{smooth structure set}\footnote{Also denoted by $\mathcal{S}^{Diff}(M)$.} $hS(M)$\label{def:smoothstructureset} of $M$ is the set of equivalence classes of simple\footnote{A homotopy equivalence is called \emph{simple} if its Whitehead torsion vanishes (see \cite[Definition 8.12]{Ra02}), which is always the case for simply connected manifolds.} homotopy equivalences $f:X^n\rightarrow M$ (sometimes called \emph{homotopy smoothings}), where $X^n$ is a smooth $n$-dimensional manifold. Two such simple homotopy equivalences $f_0:X_0\rightarrow M$ and $f_1:X_1\rightarrow M$ are \emph{equivalent} if there exists a diffeomorphism $\phi:X_0\rightarrow X_1$ such that $f_1\circ \phi \simeq f_0$.

An element $[f]\in hS(M)$ determines a normal map in the following way. Let $g:M\rightarrow X$ be a homotopy inverse of $f:X\rightarrow M$. Taking $\xi=g^*\nu_X$, we get a stable vector bundle over $M$ with a bundle map $b:\nu_X\rightarrow \xi$, and thus a normal invariant $\alpha(f)$ corresponding to $f$ (we also denote it by $\alpha(X)$ if there can be no confusion). This gives a map $\alpha:hS(M)\rightarrow [M,G/O]$ (see \cite[\textsection III.1.3.]{LdM71}).

Let $\mathcal{N}_\alpha(Q^7)$ be the restriction of $\text{Im}(\alpha)$ to Milnor projective spaces.

\begin{proposition}\label{prop:norminveellskuipmilnorquot}
The map \begin{align*}
\beta : & \text{ }  \mathcal{N}_\alpha(Q^7)\rightarrow \Z_4\\ &
\alpha(Q^7_k)\mapsto 28\mu(M^7_k)\emph{ mod }4
\end{align*} is a bijection, where we take $28\mu(M^7_k)\in\{0,1,2,...,27\}$.
\begin{proof}
Let $k_1,k_2\in\Z$. If $\alpha(Q^7_{k_1})=\alpha(Q^7_{k_2})$, then by \cite[Lemma 5.5.1. and (5.5.2)]{Ka81} and \cite[Corollary 5.4.11]{Ka81}, it follows that $(2k_1-1)\equiv \pm (2k_2-1)\text{ mod }16$. A quick computation then shows that $28\mu(M^7_{k_1})\equiv 28\mu(M^7_{k_2})\text{ mod } 4$ and thus $\beta$ is well-defined. Surjectivity of $\beta$ is immediate (take for example $k=1,2,3,4$). By the proof of \cite[V.6 Theorem]{LdM71}, the set $\mathcal{N}_\alpha(Q^7)$ has four elements. Hence the bijectivity of $\beta$ follows.
\end{proof}
\end{proposition}

Note that this result can also be deduced from the work of Mayer \cite{Ma70}.

\subsection{Eells-Kuiper invariant of the Milnor and Shimada projective spaces}

To complete the diffeomorphism classification of Milnor projective spaces, one can compute their Eells-Kuiper invariant. The computation was carried out by Tang and Zhang \cite{TZ14}, using the formula below, and the same argument can be applied to Shimada projective spaces.

Let $M^{4k-1}$ be a closed spin manifold, equipped with a Riemannian metric $g_M$ and suppose that $H^{4i}(M;\R)=0$ for all $0<i<k$. This means that there exist forms $\hat{p}_i(M)\in \Omega^{4i-1}(M)/\text{Im}(d)$\label{def:hatp} such that $p_i(M)=d\hat{p}_i(M)$ where $p_i(M)$ are the Pontrjagin forms of $M$ with respect to the metric $g_M$. Now let $\alpha(M)\in H^{4k-1}(M;\R)=\Omega^{4k-1}(M)/\text{Im}(d)$ be defined as $$\hat{A}_k(p_1,...,p_{k-1},0)-t_k L_k(p_1,...,p_{k-1},0),$$ with one factor $p_i(M)$ replaced by $\hat{p}_i(M)$ in each monomial\footnote{For example, if $\hat{A}_2(p_1(M),0)-t_2L_2(p_1(M),0)=\frac{1}{2^7\cdot 7}p_1(M)\wedge p_1(M)$ then $\alpha(M)=\frac{1}{2^7\cdot 7}p_1(M)\wedge \hat{p}_1(M)$.}. Remember that $t_k=\hat{A}_k(0,...,0,1)/L_k(0,...,0,1)$.

Then there is the following formula for the Eells-Kuiper invariant, which is helpful in the case a spin coboundary cannot be found.

\begin{theorem}\label{thm:eellskuipindexformula}
\cite[Theorem 4.8]{Goe12} Let $M$ be as above, $D_M$ denote its $Spin$ Dirac operator (see \textsection \ref{section:spindiracoperator}) and $B^{ev}_M$ its odd signature operator (see \textsection \ref{section:signatureoperator}). Then
$$\mu(M)=\frac{1}{a_k}\Big(\frac{\eta(D_M)+h(D_M)}{2}-t_k\eta(B^{ev}_M)-\int_M\alpha(M)\Big)\in \Q/\Z,$$ where $a_k=4/(3+(-1)^k)$, $\eta(D_M)$ and $\eta(B^{ev}_M)$ are the corresponding eta-invariants defined via Equation (\ref{def:etainv}) and $h(D_M)=\emph{dim}(kerD_M)$.
\end{theorem}

Note in particular that the Eells-Kuiper invariant does not depend on the choice of a Riemannian metric.

\begin{lemma}\label{lemma:eellskuipinvquotients}

The Eells-Kuiper invariant of $Q^{8n-1}_k$, $n=1,2$, is given by 
\begin{equation}\label{eqn:eellskuipinvquo}
\mu(Q^{8n-1}_k)\equiv \Big(\frac{k(k+1)}{2^{4n}\cdot q_n}\pm \frac{(2k+1)}{2^{4n+1}}\Big) \emph{ mod }1,
\end{equation} where $q_1=7$ for Milnor projective spaces and $q_2=127$ for Shimada projective spaces.

\begin{proof}

We can apply apply Theorem \ref{thm:eellskuipindexformula} to $Q^{8n-1}_k$ for both $n=1,2$. Let $D^+_W$, $D_M$ and $D_Q$ denote the corresponding Dirac operators on $W:=W^{8n}_k$, $M:=M^{8n-1}_k$ and $Q:=Q^{8n-1}_k$ from Remark \ref{rmk:spinstructmilshimproj}. If $B^{ev}_Q$ the odd signature operator of $Q$, then

$$\mu(Q)\equiv \frac{\eta(D_Q)+h(D_Q)}{2}-t_{2n} \eta(B^{ev}_Q)-c_{2n}\int_{Q}p_n(Q)\wedge \hat{p}_n(Q) \text{ mod }1$$ where $t_2=-1/(2^5\cdot 7)$, $c_2=1/(2^7\cdot 7)$, $t_4=-1/(2^9\cdot 127)$ and $c_4=1/(2^{11}\cdot 3^2\cdot 127)$. Remember that $\hat{p}_n(Q)$ is a $(4n-1)$-form satisfying $d\hat{p}_n(Q)=p_n(Q)$. Applying Theorem \ref{thm:etainvcov} to the covering $\pi:M\rightarrow Q$ with the trivial representation of $\pi_1(Q)$, we get

$$\eta(D_Q)=\frac{1}{2}\Big(\eta(D_M)+\eta_\tau(D_M)\Big),$$
$$\eta(B^{ev}_Q)=\frac{1}{2}\Big(\eta(B^{ev}_M)+\eta_\tau(B^{ev}_M)\Big),$$ where $B^{ev}_M$ is the lifted odd signature operator on $M$. Remember that the fixed point set of the action of $\tau$ on $W:=W^{8n}_k$ is the zero-section $S^{4n}$. By Theorem \ref{thm:spinapsthm}, \ref{thm:equivindthmspin} and the above, we therefore obtain

$$\eta(D_Q)=-\text{index}(D^+_W)+\int_W \hat{A}(p)-\frac{h(D_M)}{2}-\text{index}(D^+_W,\tau)+a_{spin}(S^{4n})-\frac{h_\tau(D_M)}{2}.$$ Similarly, by Theorem \ref{thm:signAPSindthm}, \ref{thm:equivindexthmsign} and the above,
$$\eta(B^{ev}_Q)=\frac{1}{2}\Big(-\text{sign}(W)+\int_W L(p)-\text{sign}(W,\tau)+a_{sign}(S^{4n})\Big).$$ Now, as we will see in \textsection \ref{sec:metricsnonnegsec}, $W$ can be equipped with a metric of non-negative scalar curvature everywhere, positive scalar curvature on the boundary $M$ and which is of product form near the boundary (see Theorem \ref{thm:metricsonmilshimprojsp} and \ref{prop:metricsons7bund}). Therefore, it follows by the vanishing theorem \ref{thm:spinvanishingthmboundary} that $\text{index}(D^+_W)$, $\text{index}(D^+_W,\tau)$, $h(D_M)$, $h_\tau(D_M)$, as well as $h(D_Q)$, all vanish.

Recall that $\text{sign}(W)=1$. Since $S^{4n}$ is the fixed point set of the action of $\tau$ on $W$, $\tau$ preserves the generator of $H^{4n}(S^{4n};\Z)$ and thus we have $\text{sign}(W,\tau)=1$.

Let $\nu_k$ be the normal bundle of the zero section $S^{4n}$ in $W$. Then $\nu_k\cong \xi_k$, where $\xi_k$ is the vector bundle associated to the $S^{4n-1}$-bundle over $S^{4n}$. By Equations (\ref{eqn:localcontrspin}) and (\ref{eqn:ahatpiofE}), we have
$$a_{spin}(S^{4})=\pm\frac{1}{2^{5}}\int_{S^{4n}}p_1(\xi_k), \qquad a_{spin}(S^{8})=\pm\int_{S^{8}}\Big(\frac{5}{2^{11}\cdot 3} p^2_1(\xi_k)-\frac{1}{2^9\cdot 3}p_2(\xi_k)\Big),$$
so that, using Theorem \ref{thm:spherebundleproperties}.4, we obtain
$$a_{spin}(S^{4n})=\pm \frac{(2k+1)}{2^{4n}}.$$

 Similarly, by Equation (\ref{eqn:localcontrsign}) and Theorem \ref{thm:spherebundleproperties}.1,
$$a_{sign}(S^{4n})=\int_{S^{4n}}e(\xi_k)=1.$$ We also have 
$$\int_{Q}p_n(Q)\wedge \hat{p}_n(Q)=\frac{1}{2}\int_{M}p_n(M)\wedge \hat{p}_n(M)$$ since $\int_M\pi^*(\omega)=deg(\pi)\int_Q\omega$ for any $(8n-1)$-form $\omega$. Finally, by \cite[Lemma 2.7]{KS93},
$$\int_W p_n(W)\wedge p_n(W)-\int_M p_n(M)\wedge\hat{p}_n(M) = \langle \overline{p}^2_n(W),[W,M]\rangle.$$ Applying Equation (\ref{eqn:pontrdiskbund}) and putting all of the above together, the result now follows.
\end{proof}

\end{lemma}

\begin{remark}
Observe that the sign in Equation (\ref{eqn:eellskuipinvquo}) depends on the choice of the $Spin$ structure on $Q^{8n-1}_k$ (see \cite[p.58]{Mi65}). Indeed, we have $H^1(Q^{8n-1}_k;\Z_2)\cong \Z_2$ and so there are two different $Spin$ structures on $Q^{8n-1}_k$. The Eells-Kuiper invariant $\mu(Q^{8n-1}_k)$ therefore has to be interpreted as a pair of values, not as a singular value.
\end{remark}



\begin{proposition}\label{prop:eellskuiperequal}
Let $Q^{7}_{k_i}=M^{7}_{k_i}/\tau$ for $i=0,1$. Then $\mu(M^{7}_{k_0})=\mu(M^{7}_{k_1})$ implies $\mu(Q^{7}_{k_0})=\mu(Q^{7}_{k_1})$.

\begin{proof}
Observe that $\mu(M^{7}_{k})=\mu(M^{7}_{k+56m})$ and $\mu(Q^{7}_{k})=\mu(Q^{7}_{k+56m})$ for $m\in\Z$. The result now follows by computing and comparing the different Eells-Kuiper invariants for $k=0,1,...,55$.
\end{proof}
\end{proposition}

\begin{proposition}\label{prop:eellskuipershimadaquodiffvalues}
There are 4096 different pairs of values for $\mu(Q^{15}_k)$.
\begin{proof}
This is achieved through use of the C++ code in the Appendix B.
\end{proof}
\end{proposition}

\subsection{Classification of Milnor projective spaces}

\begin{theorem}\label{thm:diffmilsphimpliesdiffquot}
Let  $Q^7_{k_i}=M^7_{k_i}/\tau$ be a Milnor projective space for $k_i\in\Z$, $i=0,1$. If $M^7_{k_0}$ is diffeomorphic to $M^7_{k_1}$, then $Q^7_{k_0}$ is diffeomorphic to $Q^7_{k_1}$.
\begin{proof}

If $M^7_{k_0}$ is diffeomorphic to $M^7_{k_1}$, then $28\mu(M^7_{k_0})=28\mu(M^7_{k_1})$ and therefore by Proposition \ref{prop:norminveellskuipmilnorquot} their normal invariants are equal: $\alpha_7(Q^7_{k_0})=\alpha_7(Q^7_{k_1})$. By Theorem \ref{thm:blhirschmilnor}, the Browder-Livesay invariant is $\sigma(M^7_k,\tau)=0$ for all $k\in\Z$, hence by Theorem \ref{thm:projspbrownorm} it follows that $Q^7_{k_0}$ is diffeomorphic to $Q^7_{k_1}\#\Sigma^7$ for some sphere $\Sigma^7\in bP_8$. Suppose $\mu(\Sigma^7)=\frac{l}{28} \text{ mod }1$, where $0\leq l<28$ is an integer. By the properties of the Eells-Kuiper invariant (see Proposition \ref{prop:eellskuipproperties}), we have

$$\mu(Q^7_{k_0})=\mu(Q^7_{k_1}\#\Sigma^7)\equiv\mu(Q^7_{k_1})+\frac{l}{28} \text{ mod } 1.$$ Proposition \ref{prop:eellskuiperequal} now implies that  $l=0$. Therefore $\Sigma^7\cong S^7$ and finally $Q^7_{k_1}\# \Sigma^7 \cong Q^7_{k_1}$.

\end{proof}
\end{theorem}

\begin{theorem}\label{thm:classmilprojsp}
There are 16 different (oriented) diffeomorphism classes of Milnor projective spaces. All of the 16 diffeomorphism types can be realized by an infinite family of such quotients.

\begin{proof}
The first statement follows from Corollary \ref{cor:numberofdifftypesspherebund}. The second statement follows from Theorem \ref{thm:diffeoclassspherebun}, \ref{thm:diffmilsphimpliesdiffquot} and Corollary \ref{cor:infinitefamilysphbund}.
\end{proof}

\end{theorem}

This completes the proof of Theorem \ref{mainthm:diffclassmilproj}.

\subsection{Diffeomorphism finiteness of Shimada projective spaces}

As of the time of writing, the normal invariants of Shimada projective spaces are still unknown. Hence, the best we can do is to give a finiteness result and a lower limit for the number of diffeomorphism types of Shimada projective spaces.

\begin{lemma}\label{lemma:diffinrp15}
There are only finitely many different oriented diffeomorphism types of Shimada projective spaces.

\begin{proof}
If $a_i$ is the order of $\pi_i(G/O)$ and $b_i$ is the order of $\pi_i(G/O)\otimes \Z_2$, then it can be shown that the order of $[\R\text{P}^{15},G/O]$ is less than or equal to $\prod_{i=1}^{14}a_{15}b_i$ (see \cite[V.1]{LdM71}). Now, since $a_{15}=2$ and $b_i$ is finite for all $i$ (see \cite{Sul96}), it follows that the order of $[\R\text{P}^{15},G/O]$ is finite. In particular, there are only finitely many distinct normal invariants that Shimada projective spaces can have.

By Theorem \ref{thm:blhirschmilnor} the Browder-Livesay invariant of a Shimada projective space $Q^{15}_k$ vanishes. Therefore, by the above and Theorem \ref{thm:projspbrownorm}, the result follows since $|bP_{16}|=8128$ is finite.
\end{proof}
\end{lemma}

\begin{proposition}\label{prop:difffinshimadaproj}
There are finitely many, but at least 4096 different oriented diffeomorphism types of Shimada projective spheres which can be realized by an infinite family of orientation preservingly diffeomorphic manifolds.
\begin{proof}
The first statement follows from Lemma \ref{lemma:diffinrp15}. The second statement follows from Proposition \ref{prop:eellskuipershimadaquodiffvalues}. The last statement now follows by considering $\{Q^{15}_{k+130048m}\}_{m\in \Z}$, which all have the same Eells-Kuiper invariant as $Q^{15}_{k}$ for any $k\in \Z$ (see Lemma \ref{lemma:eellskuipinvquotients}).
\end{proof}
 
\end{proposition}

\section{Construction of the Riemannian metrics}

\subsection{Metrics of non-negative sectional curvature}\label{sec:metricsnonnegsec}

Consider princiapl $S^3\times S^3$-bundles over $S^4$. These bundles are classified by elements of $\pi_3(SO(4))=\Z\oplus \Z$. Let $P^{10}_{k,l}$\label{def:P10cohomogenonemanifold} denote the total space corresponding to $[\alpha]\in \pi_3(SO(4))$ where\footnote{Note that Grove and Ziller \cite{GZ00} use different indices.} $\alpha(q)u=q^{k+l} u q^{-k}$.

The manifold $P^{10}_{k,l}$ admits a cohomogeneity one action by $S^3\times S^3\times S^3/\pm(1,1,1)$ with codimension 2 singular orbits \cite[Proposition 3.11]{GZ00}. By \cite[Theorem E]{GZ00}, it therefore admits an $S^3\times S^3\times S^3$-invariant metric with non-negative sectional curvature.

Let $S^3\times S^3\times 1$ be the subaction of $S^3\times S^3\times S^3$ corresponding to the principal $S^3\times S^3$-bundle action on $P^{10}_{k,l}$. Note that $S^3\times S^3\times 1$ acts freely and isometrically on $P^{10}_{k,l}$ and that the quotient $P^*_{k,l}:=P^{10}_{k,l}/(-1,-1)$ is the total space of the associated principal $SO(4)$-bundle over $S^4$. Let $S^3\times S^3$ act on $S^3$ via $(q_1,q_2)\cdot v=q_1vq^{-1}_2$, where quaternion multiplication is understood. Then $M^7_{k,l}:=P^{10}_{k,l}\times_{S^3\times S^3} S^3$ is an $S^3$-bundle over $S^4$ with structure group $SO(4)$ and Euler class $l$.

\begin{theorem}\label{thm:metricsonmilshimprojsp}
Fix $k\in\Z$, $l>0$ and let $W^8_{k,l}$ be the disk bundle associated to $M^7_{k,l}$. Then $M^7_{k,l}$ admits a metric that is simultaneously of non-negative sectional curvature and positive scalar curvature, which will be called its \emph{Grove-Ziller metric} and denoted by $\tilde{g}^{GZ}_{k,l}$\label{def:grovezillermilnorsphere}. It extends to a metric $h_{k,l}$ of non-negative sectional curvature on $W^8_{k,l}$ which is of product form near the boundary.

When $l=1$, the metric $\tilde{g}^{GZ}_k:=\tilde{g}^{GZ}_{k,1}$ on the Milnor sphere $M^7_k$ descends to a metric of non-negative sectional curvature on $Q^7_k$, which we likewise call its \emph{Grovez-Ziller metric} and denote by $g^{GZ}_k$\label{def:grovezillermilnorproj}. It satisfies $\tilde{g}^{GZ}_k=\pi_k^*(g^{GZ}_k)$, where $\pi_k:M^7_k\rightarrow Q^7_k$ is the projection.

\begin{proof}

The construction of the metrics on $M^7_{k,l}$ and $W^8_{k,l}$ with the above properties has been discussed in \cite{De17} (see also \cite[Theorem 6.1.2]{We21}).

The involution $\tau$ which defines the Milnor projective spaces is induced by the action of $Id_{P}\times (-Id_{S^3})$ on $P^{10}_{k,1}\times S^3$. Therefore, the metric $\tilde{g}^{GZ}_{k}$ is invariant under $\tau$ and $Q^7_k:=M^7_k/\tau$ inherits a metric $g^{GZ}_k$ with non-negative sectional curvature, satisfying the required properties.

\end{proof}
\end{theorem}



\subsection{Positive Ricci curvature metrics}\label{sec:posriccicurv}

We can now construct positive Ricci curvature metrics on $S^{4n-1}$-bundles over $S^{4n}$ that extend to the corresponding disk bundle.

\begin{proposition}\label{prop:metricsons7bund}

Let $M^{8n-1}_{k,l}$ be the total space of a linear $S^{4n-1}$-bundle over $S^{4n}$ and $W^{8n}_{k,l}$ the total space of the associated disk bundle, $n=1,2$. There exists a metric $\tilde{g}_{k,l}$\label{def:metricricciposdiskbundle} on $W^{8n}_{k,l}$ which has positive scalar curvature, is of product form near the boundary $M^{8n-1}_{k,l}$ and such that $g_{k,l}=\tilde{g}_{k,l}|_{M^{8n-1}_{k,l}}$\label{def:metricricciposspherebundle} has positive Ricci curvature.
\begin{proof}

Let $P:=P_{k,l}$, $M:=M^{8n-1}_{k,l}$ and $W:=W^{8n}_{k,l}$ be the total space of the associated principal $SO(4n)$-bundle, the $S^{4n-1}$-bundle and the $D^{4n}$-bundle over $S^{4n}$ respectively. Equip $S^{4n}$ with the round metric $g_R$ and $D^{4n}$ with a torpedo metric $g_{tor}$ (see \cite{Wa11} for the definition and properties).

Fix a connection on $P\rightarrow S^{4n}$. Then by Vilms \cite{Vi70} (see also \cite[Proposition 2.7.1]{GW09}), there exists a unique metric on $W=P\times_{SO(4n)} D^{4n}$ such that the projection $\pi_W:W\rightarrow S^{4n}$ is a Riemannian submersion with totally geodesic fibers isometric to $(D^{4n},g_{tor})$. Denote this metric by $\tilde{g}$. The restriction $g:=\tilde{g}|_M$ on $M=P\times_{SO(4n)}S^{4n-1}$ corresponds to metric of Vilms applied to $M=P\times_{SO(4n)}S^{4n-1}$.

Let $\mathcal{V}:=\text{ker}(\pi_W)_*$ be the vertical distribution on $W$ and set $\mathcal{H}$ as the orthogonal complement of $\mathcal{V}$ with respect to $\tilde{g}_W$. Consider the canonical variation $$\tilde{g}_t|_{\mathcal{V}}:=t\cdot \tilde{g}|_{\mathcal{V}}, \qquad \tilde{g}_t|_{\mathcal{H}}:=\tilde{g}|_{\mathcal{H}} \qquad \text{and} \qquad \tilde{g}_t(\mathcal{V},\mathcal{H}):=0$$ on $W$, where $t\in \R_{\geq 0}$. This amounts to setting 
\begin{equation}\label{eqn:canmetric}
\tilde{g}_t(X,Y):=t\cdot g_{tor}(X^\mathcal{V},Y^\mathcal{V})+\pi_W^*g_{R}(X,Y), \qquad \text{for }X,Y\in T_xW.
\end{equation}
If we set $g_t:=(\tilde{g}_t)|_M$ and restrict to $\mathcal{V}_M:=\text{ker}(\pi_M)_*\subset \mathcal{V}$ with its corresponding horizontal distribution $\mathcal{H}_M$, then this simultaneously corresponds to a canonical variation on $M=\partial W $:
$$g_t|_{\mathcal{V}_M}=t\cdot g|_{\mathcal{V}_M}, \qquad g_t|_{\mathcal{H}_M}=g|_{\mathcal{H}_M} \qquad \text{and} \qquad g_t(\mathcal{V}_M,\mathcal{H}_M)=0.$$
Now by \cite[Theorem 2.7.3]{GW09} and its proof applied to $E=M$, $B=S^{4n}$ and $F=S^{4n-1}$, there is an $0<\epsilon<<1$ such that $g_M:=g_\epsilon$ is of positive Ricci curvature.

Next we show that the metric on $W$ has positive scalar curvature. By \cite[9.70(d)]{Be87}, the scalar curvature of the canonical variation metric is given by $$scal_{\tilde{g}_t}=\frac{1}{t}scal_{g_F}+scal_{g_B}\circ \pi_W - t|A|^2$$ where $A$ is a tensor field on $W$ and in our case, $g_F=g_{tor}$ and $g_B=g_R$. Obviously, $scal_{g_{tor}}>0$ and $scal_{g_R}>0$. Therefore, choosing $\epsilon$ to be even smaller if necessary, the metric $g_W:=g_\epsilon$ is of positive scalar curvature everywhere and restricts to the positive Ricci curvature metric $g_M=g_W|_M$. 

Finally, Equation (\ref{eqn:canmetric}) shows that $g_W$ is of product form near the boundary. As we have mentioned above, the fibers on $W$ are isometric to $(D^{4n},g_{tor})$ and the canonical variation corresponds to shrinking the fibers, therefore respecting the product form near the boundary (see \cite[p.10]{Kor20}).

\end{proof}
\end{proposition}

\subsection{Metrics on Shimada projective spaces}\label{sec:metricshimadaproj}

Let $M^{15}_k$ be a Shimada sphere equipped with the metric $g_k$ from Proposition \ref{prop:metricsons7bund} (remember that $M^{15}_k:=M^{15}_{k,1}$). Since the fibers of the Riemannian submersion $\pi_S:M^{15}_k\rightarrow S^8$ are isometric to the round sphere $(S^7,g_R)$ (see the proof of Proposition \ref{prop:metricsons7bund}), it follows that the involution $\tau$ on $M^{15}_k$, which is induced by fiberwise antipodal maps, is an isometry. Therefore, the induced metric $g'_k$ on the quotient $Q^{15}_k:=M^{15}_k/\tau$ is of positive Ricci curvature and satisfies $g_k=\pi_k^*(g'_k)$\label{def:metricricciposshimadaproj} where $\pi_k:M^{15}_k\rightarrow Q^{15}_k$ is the canonical projection.

\begin{remark}\label{rem:posriccimiln}
The exact same argument applies to construct positive Ricci curvature metrics on Milnor projective spaces.
\end{remark}

\section{Moduli spaces of Riemannian metrics and proofs}\label{sec:modspaceproof}

For more details on moduli spaces of Riemannian metrics, see \cite{TW15}.

Let $M$ be a compact smooth manifold and $\mathcal{R}(M)$\label{def:setofriemanmetric} the set of Riemannian metrics on $M$, equipped with the $\mathcal{C}^\infty$-topology of uniform convergence of all the derivatives. If we restrict to metrics with $sec\geq 0$, $Ric>0$ and $scal>0$, we get the corresponding sets $\mathcal{R}_{sec\geq 0}(M)$, $\mathcal{R}_{Ric>0}(M)$ and $\mathcal{R}_{scal>0}(M)$. The group of diffeomorphisms Diff($M$)\label{def:groupofdiffeom} acts on $\mathcal{R}(M)$ by taking pullbacks of the metrics. The \textit{moduli space of Riemannian metrics} of $M$ is defined as the quotient space $\mathcal{M}(M):=\mathcal{R}(M)/\text{Diff}(M)$\label{def:modulispace} whose elements are isometry classes of Riemannian metrics. If we restrict to isometry classes of metrics of $sec\geq 0$, $Ric>0$ and $scal>0$ we get corresponding moduli spaces $\mathcal{M}_{sec\geq 0}(M)$, $\mathcal{M}_{Ric> 0}(M)$ and $\mathcal{M}_{scal> 0}(M)$.

Let $C$ either be $sec\geq 0$ or $Ric>0$.

\begin{proposition}\label{prop:prooftool}

Let $M^{4k-1}$ be a closed spin $(4k-1)$-dimensional manifold. Let $\phi_i:M\rightarrow M_i$ be an orientation preserving diffeomorphism and $W_i^{4k}$ a spin manifold with boundary $\partial W_i = M_i$ for $i=0,1$. Assume that all the spin structures are unique. Let $\tilde{g}_i$ be a metric on $W_i$ which is of $scal\geq 0$ everywhere, $scal>0$ on the boundary, of product form near the boundary and such that $g_i:=\tilde{g}_i|_{M_i}$ satisfies the curvature condition $C$. Let $h_i:=\phi^*_i(g_i)$, $i=0,1$. Then if $h_0$ and $\psi^*(h_1)$ lie in the same path component of $\mathcal{R}_C(M)$ for some orientation preserving $\psi\in \emph{Diff}(M)$, there exist $a>>0$ and a path $\gamma:I\rightarrow \mathcal{R}_{scal>0}(M)$ with $\gamma(0)=h_0$ and $\gamma(1)=\psi^*(h_1)$ such that the spin manifold

$$(X^{4k},g):=(W_0,\tilde{g}_0)\cup_{\phi^{-1}_{0}}(M\times[0,a],\gamma(t/a)+dt^2)\cup_{\phi_{1}\psi}(-W_1,\tilde{g}_1)$$ satisfies $\emph{index}(D^+_X)=0$ and $\emph{sign}(X)=0$, where $D^+_X$ is the $Spin^+$ Dirac operator of $X$.

\begin{proof}
Let $\gamma:I\rightarrow \mathcal{R}_C(M)$ be the path with endpoints $\gamma(0)=h_0$ and $\gamma(1)=\psi^*(h_1)$.

If $C=\{Ric>0\}$, the path $\gamma$ in particular also lies in $\mathcal{R}_{scal> 0}(M)$. If $C=\{sec\geq 0\}$, by Böhm and Wilking \cite{BW07}, the path $\gamma$ instantly evolves to a path in $\mathcal{R}_{Ric>0}(M)$ under the Ricci flow. If we concatenate this path with the orbits of the endpoints of $\gamma$ under the Ricci flow, we get a path in $\mathcal{R}_{scal> 0}(M)$ with the same endpoints as $\gamma$. We will still denote this resulting path by $\gamma$. In any case, we can now reparametrize this path (and still denote it by $\gamma$) in such a way that it becomes constant near the endpoints $\gamma(0)=\phi^*_0(g_0)$ and $\gamma(1)=\psi^*(\phi^*_1(g_1))$.

According to Gromov and Lawson \cite[Lemma 3]{GL80}, the product $M\times [0,a]$ equipped with $\gamma(t/a)+dt^2$ has positive scalar curvature for some $a\gg0$. Hence, we can define $X$ as above. The index of the $Spin^+$ Dirac operator of $X$ corresponds to the relative index invariant of Gromov and Lawson \cite[p.116]{GL83}.

Since the metric on $X$ has $scal\geq 0$ everywhere and $scal>0$ on the cylinder, by the standard argument of Lichnerowicz \cite{Li63}, the index of $D^+_X$ vanishes. Note that $X$ is diffeomorphic to $W_0\cup_{\alpha}(-W_1),$ where $\alpha:=\phi_{1}\psi\phi^{-1}_{0}$. Then, by a formula which is due to Novikov (see \cite[Proposition (7.1)]{ASIII68}), we have $\text{sign}(X)=\text{sign}(W_0)+\text{sign}(-W_1)=0$.

\end{proof}

\end{proposition}

\subsection{Proof of Theorems \ref{mainthm:s3bundles} and \ref{mainthm:s7bundles}}

Fix $k\in\Z$, $l>0$ and let $M^{8n-1}:=M^{8n-1}_{k,l}$ be the total space of a linear $S^7$-bundle over $S^8$. Then, by Corollary \ref{cor:infinitefamilysphbund}, for $k(m)=k+2^{4n-2}lm\cdot q_n$, $m\in\Z$, the set $\{M^{8n-1}_{k(m),l}\}_{m\in \Z}$ is an infinite family of manifolds orientation preservingly diffeomorphic to $M^{8n-1}$.

Fix $m_0,m_1\in \Z$ such that $|2k_0+l|\neq |2k_1+l|$, where $k_i:=k+2^{4n-2}lm_i\cdot g_n$ for $i=0,1$. Denote by $\phi_{i}:M\rightarrow M^{8n-1}_{k_i,l}$ the diffeomorphism for $i=0,1$ (see Theorem \ref{thm:diffeoclassspherebun}). For $n=1$, equip $M^{7}_i:=M^7_{k_i,l}$ with the metric $g_{i}:=\tilde{g}^{GZ}_{k_i,l}$ from Theorem \ref{sec:metricsnonnegsec} and if $n=2$, equip $M^{15}_i:=M^{15}_{k_i,l}$ with the metric $g_{i}:=g_{k_i,l}$ from Proposition \ref{prop:metricsons7bund}. In both cases, let $h_i:=\phi_i^*(g_i)$ for $i=0,1$. Denote $W^{8n}_i:=W^{8n}_{k_i,l}$.

The proof goes by contradiction. Assume there is a path $\delta$ in $\mathcal{M}_{C}(M^{8n-1})$ (where $C=\{sec\geq 0\}$ if $n=1$ and $C=\{Ric>0\}$ if $n=2$) with endpoints $\delta(0)=[h_0]$ and $\delta(1)=[h_1]$. As a consequence of the Ebin slice theorem (see \cite{Eb70} and \cite[Proposition 4.6]{CK19}), this path lifts to a path $\epsilon$ in $\mathcal{R}_{C}(M^{8n-1})$ such that $\epsilon(0)=\phi^*_0(g_0)$ and $\epsilon(1)=\psi^*(\phi^*_1(g_1))$ for some $\psi\in\text{Diff}(M^{8n-1})$. If $\psi$ is orientation reversing, we can replace $g_{1}$ by its pullback under an orientation reversing diffeomorphism of $M^{8n-1}_1$ (the pullback of $g_1$ by this orientation reversing diffeomorphism still is a representative of $[\phi^*_1(g_1)]$ in $\mathcal{M}_{C}(M^{8n-1})$), in order to compensate. Hence, we can from now on assume (without loss of generality) that $\psi$ is orientation preserving.

Then, by Proposition \ref{prop:prooftool}, we have $\text{index}(D^+_X)=0$ and $\text{sign}(X)=0$, where $X^{8n}:=W_0\cup_{\phi^{-1}_{0}}(M\times[0,a])\cup_{\phi_{1}\psi}(-W_1)$.

By the Atiyah-Singer index theorem \cite[Theorem (5.3)]{ASIII68}, we have $$\text{index}(D^+_X)=\langle\hat{A}(X^8),[X^8]\rangle=\langle\frac{-4p_2(X^8)+7p^2_1(X^8)}{5760},[X^8]\rangle$$ and
\begin{align}\label{eqn:indexX}
\text{index}(D^+_X)&=\langle\hat{A}(X^{16}),[X^{16}]\rangle \\
&=\langle \frac{-192p_4+512p_3p_1+208p^2_2-904p_2p^2_1-904p^4_1}{464486400},[X^{16}]\rangle.
\end{align} By Hirzebruch's signature theorem  (see \cite[Theorem (6.6)]{ASIII68}), $$\text{sign}(X^8)=\langle L(X^8),[X^8]\rangle =\langle \frac{7p_2(X^8)-p^2_1(X^8)}{45},[X^8]\rangle$$ and \begin{equation}\label{eqn:signX}
\text{sign}(X^{16})=\langle L(X^{16}),[X^{16}]\rangle=\langle \frac{381p_4-71p_3p_1-19p^2_2+22p_2p^2_1-3p_1^4}{14175},[X^{16}]\rangle.
\end{equation}

If $n=1$, the two preceding constraints imply $\langle p^2_1(X^8),[X^8]\rangle=0$. If $n=2$, by the Mayer-Vietoris exact sequence, we have $H^4(X^{16};\Z)=0$. Therefore $p_1(X^{16})=0$ and so by the above, Equations (\ref{eqn:indexX}) and (\ref{eqn:signX}) reduce to
$$\langle -192p_4+208p^2_2,[X^{16}]\rangle=0,$$
$$\langle 381p_4-19p^2_2,[X^{16}]\rangle=0.$$ It follows that both $\langle p_4,[X^{16}]\rangle$ and $\langle p^2_2,[X^{16}]\rangle$ must vanish.

However, using Equation (\ref{eqn:pontrdiskbund}), for we compute $$\langle p^2_1(X^{8n}),[X^{8n}]\rangle=\langle(\overline{p}^2_1(W^{8n}_{0}),[W^{8n}_{0},M^{8n-1}_{0}]\rangle-\langle(\overline{p}^2_1(W^{8n}_{1}),[W^{8n}_{1},M^{8n-1}_{1}]\rangle$$ $$ = \frac{(4n-2)^2}{l}\Big((2k_0+l)^2-(2k_1+l)^2\Big).$$ See \textsection \ref{sec:diffclasssphbund} for the definition of $\overline{p}_1(W_{i})$ for $i=0,1$.

This is a contradiction, since we assumed $|2k_0+l|\neq |2k_1+l|$ at the beginning of the proof. Hence $[h_0]$ and $[h_1]$ cannot lie in the same path component of $\mathcal{M}_{C}(M^{8n-1})$.

In case $n=1$, one can use use the exact same argument but starting with the positive Ricci curvature metrics from Proposition \ref{prop:metricsons7bund} to arrive at the same conclusion on $\mathcal{M}_{Ric>0}(M^{7})$.

Using Theorem \ref{thm:diffeoclassspherebun}, this completes the proof of Theorems \ref{mainthm:s3bundles} and \ref{mainthm:s7bundles}. \proofend 

\begin{remark}
If $l=1$, i.e $M^{15}_{k,1}$ is a homotopy 15-sphere, this result was already proved by Wraith \cite{Wr11} using a different method to construct suitable positive Ricci curvature metrics that extend to a coboundary (see also \cite{Wr97}). Wraith's method can also be applied to the moduli space of more general $S^7$-bundles over $S^8$, leading to the same result.
\end{remark}

\subsection{Proof of Theorems \ref{mainthm:milnprojsp} and \ref{mainthm:shimprojsp}}

Fix $k\in\Z$ and let $Q^{8n-1}:=Q^{8n-1}_k=M^{8n-1}_k/\tau$ be a Milnor projective space for $n=1$ or a Shimada projective space for $n=2$ and set $M^{8n-1}:=M^{8n-1}_k$.

Fix $m_0, m_1\in \Z$ such that $|2k_0+1|\neq |2k_1+1|$ (where $k_i=k+2^{4n-2}m_i\cdot q_n$ for $i=0,1$) and such that there exist orientation preserving diffeomorphisms $\Psi_i:Q^{8n-1}\rightarrow Q^{8n-1}_{k_i}$ and $\tau$-equivariant diffeomorphisms $\tilde{\Psi}_i:M^{8n-1}\rightarrow M^{8n-1}_{k_i}$ for $i=0,1$.

Equip $Q_i^7$ and the corresponding Milnor sphere $M_i^7$ with the above Grove-Ziller metric $g_i:=g^{GZ}_{k_i}$ and $\tilde{g}_i:=\tilde{g}^{GZ}_{k_i}$ respectively (see Theorem \ref{thm:metricsonmilshimprojsp}). Equip $Q^{15}$ and the corresponding Shimada sphere $M^{15}:=M^{15}_k$ with the metric $g_i:=g'_{k_i}$ from section \ref{sec:metricshimadaproj} and $\tilde{g}_i:=g_{k_i,1}$ from Proposition \ref{prop:metricsons7bund} respectively.

We then have the following commutative diagram.

\begin{center}
\begin{tikzcd}
  M^{8n-1}_{0}  \arrow[d, "\pi_{0}" '] & \arrow[l, "\tilde{\Psi}_{0}" '] M^{8n-1} \arrow[r, "\tilde{\Psi}_{1}"] \arrow[d, "\pi"]  & M^{8n-1}_{1} \arrow[d, "\pi_{1}"] \\ 
  Q^{8n-1}_{0} &  \arrow[l, "\Psi_{0}" '] Q^{8n-1} \arrow[r, "\Psi_{1}" ]  & Q^{8n-1}_{1}
  \end{tikzcd}
\end{center} where $\pi_i$ and $\pi$ are the corresponding canonical projections for $i=0,1$.

The proof is by contradiction. Let $h_i:=\Psi_{i}^*(g_{i})$ for $i=0,1$. Consider $[h_0],[h_1]\in\mathcal{M}_{C}(Q^{8n-1})$ (where $C=\{sec\geq 0\}$ if $n=1$ and $C=\{Ric>0\}$ if $n=2$) and assume that there is a path between them in this moduli space. As a consequence of the Ebin slice theorem, this path lifts to a path $\epsilon$ in $\mathcal{R}_{C}(Q^{8n-1})$ such that $\epsilon(0)=h_0$ and $\epsilon(1)=\phi^*(h_1)$ for some $\phi\in\text{Diff}(Q^{8n-1})$. If $\phi$ is orientation reversing, we can replace $g_1$ by its pullback under another orientation reversing diffeomorphism of $Q_1^{8n-1}$, in order for the composition of the two diffeomorphisms to be orientation preserving. Hence, we can from now on assume (without loss of generality) that $\phi$ is orientation preserving.

Recall that $\pi$ is a local isometry. Then, the pullback $\tilde{\epsilon}:=\pi^*(\epsilon)$ is a path in $\mathcal{R}_{C}(M^{8n-1})$ starting at $\tilde{\epsilon}(0)=\pi^*(h_0)$ and ending at $\tilde{\epsilon}(1)=\pi^*\phi^*(h_1)$. Now as a special case of \cite[(3.2) Proposition.]{Bro67}, there exists an orientation preserving diffeomorphism $\tilde{\phi}\in\text{Diff}(M^7)$ such that $\pi\circ \tilde{\phi}=\phi\circ \pi$. By commutativity of the above diagrams, we see that we can rewrite the endpoints as $\tilde{\epsilon}(0)=\tilde{\Psi}_{0}^*\pi_{0}^*(g_{0})=\tilde{h}_0$ and $\tilde{\epsilon}(1)=\tilde{\phi}^*\tilde{\Psi}_{1}^*\pi_{1}^*(g_{1})=\tilde{\phi}^*(\tilde{h}_1)$ where $\tilde{h}_i:=\tilde{\Psi}_{i}^*(\tilde{g}_{i})$ for $i=0,1$.

We can now apply Proposition \ref{prop:prooftool} to conclude that $\text{index}(D^+_X)=0$ and $\text{sign}(X)=0$, where $X^{8n}:=W_0\cup_{\phi^{-1}_{0}}(M\times[0,a])\cup_{\phi_{1}\psi}(-W_1)$.

The considerations and computations of the proof of Theorems \ref{mainthm:s3bundles} and \ref{mainthm:s7bundles} apply in exactly the same way (remember that $l=1$) in this situation, and hence we arrive at the desired contradiction. Therefore $[h_0]$ and $[h_1]$ cannot lie in the same path component of $\mathcal{M}_{C}(Q^{8n-1})$.

By Theorem \ref{thm:classmilprojsp}, for the 16 different diffeomorphism types there are infinitely many indices such that the Milnor projective spaces are all pairwise diffeomorphic, hence it follows that the moduli space of $sec\geq 0$ metrics has infinitely many path components and thus Theorem \ref{mainthm:milnprojsp} has been proved. 

By Proposition \ref{prop:difffinshimadaproj}, there are at least 4096 different diffeomorphism types which can be expressed as an infinite family of diffeomorphic manifolds with different indices. Hence, by the above, for each of these diffeomorphism types the moduli space of $Ric>0$ metrics has infinitely many path components and thus Theorem \ref{mainthm:shimprojsp} follows. \proofend


\section{Appendix A. Index theory}

Let $W$ be a compact Riemannian manifold which is of product form near its boundary $M=\partial W$. Let $\mathcal{D}_W:\Gamma(E)\rightarrow \Gamma(F)$ be a linear first order elliptic differential operator for some vector bundles $E$ and $F$ over $W$, endowed with a smooth inner product. Denote by $\mathcal{D}^*_W$ its adjoint operator.
Let $G$ be a subgroup of the isometry group of $W$ and assume that the action of $G$ is a product near the boundary. Suppose furthermore that the action lifts to $E$ and $F$ and that the induced map on sections commutes with $\mathcal{D}_W$. Then, for each $g\in G$, the equivariant index of $\mathcal{D}_W$ is defined by 
\begin{equation}\label{def:equivindex}
\text{index}(\mathcal{D}_W,g):=tr(g|_{ker \mathcal{D}_W})-tr(g|_{ker \mathcal{D}^*_W}).
\end{equation} We also set
\begin{equation}\label{def:equivkernel}
h_g(\mathcal{D}_M):=tr(g|_{ker \mathcal{D}_M}).
\end{equation}

The \emph{equivariant eta-invariant} $\eta_g(\mathcal{D}_M)$ is defined as $\eta_g(0)$, where \begin{equation}\label{eqn:equivetainv}
\eta_g(z):=\sum_{\lambda\neq 0} \frac{\text{sign}(\lambda) tr(g^\#_\lambda)}{|\lambda|^z},
\end{equation} for $g^\#_\lambda$ the map induced by $g$ on the eigenspaces $E'_\lambda$ of $\mathcal{D}_M$.

For $g=e=Id_W$, we get
\begin{equation}\label{def:index}
\text{index}(\mathcal{D}_W):=\text{index}(\mathcal{D}_W,e)=\text{dim}(ker\mathcal{D}_W)-\text{dim}(ker\mathcal{D}^*_W),
\end{equation}
\vspace{-3mm}
\begin{equation}\label{def:kernel}
h(\mathcal{D}_M):=h_{e}(\mathcal{D}_M)=\text{dim}(ker\mathcal{D}_M),
\end{equation} and
\vspace{-2mm}
\begin{equation}\label{def:etainv}
\eta(\mathcal{D}_M):=\eta_e(\mathcal{D}_M)=\eta(0)
\end{equation} is the \emph{eta-invariant} associated to $\mathcal{D}_M$.

\begin{remark}\label{rmk:normalbundledecomp}
Let $W$ be a compact, oriented Riemannian manifold (possibly with boundary) and $g:W\rightarrow W$ an orientation preserving isometry. Consider a fixed point component $N\subset W^g$ and denote by $\nu$ its normal bundle in $W$. Then the differential of $g$ induces a bundle isometry $dg:\nu\rightarrow \nu$ and from representation theory, it follows immediately that there is a direct sum decomposition \label{def:normalbunddecomposition} $$\nu=\nu(\pi)\oplus \bigoplus_{0<\theta <\pi}\nu(\theta)$$ where $\nu(\pi)$ is real and $\nu(\theta)$ is complex for $0<\theta <\pi$, $dg$ acts on $\nu(\pi)$ via multiplication by $-1$ and on $\nu(\theta)$ via multiplication by $e^{i\theta}$ (see for example \cite[pp.262-265]{LM89}).
\end{remark}

\subsection{$Spin$ Dirac operator}\label{section:spindiracoperator}

See \cite[Chapter II]{LM89} and \cite[Chapter 11]{Ni07} for more details on the $Spin$ Dirac operator.

Let $M^{n}$ be an $n$-dimensional compact spin manifold (possibly with non-empty boundary). Let $S$ be the associated spinor bundle. Equip $S$ with the Riemannian connection $\nabla$ induced by the canonical Riemannian connection on $P_{SO}$. The \emph{$Spin$ Dirac operator}\label{def:spindiracop} of $S$ at $x\in M$ is the first-order differential operator $D_M:\Gamma(S)\rightarrow \Gamma(S)$ defined by $D_M(\sigma):=\sum_{j=1}^n e_j\cdot \nabla_{e_j}\sigma$ where $\{e_1,...,e_n\}$ is an orthonormal basis of $T_xM$ and ``$\cdot$" denotes Clifford multiplication. It is well-known that this operator is elliptic and formally self-adjoint (see \cite[II\textsection 5]{LM89}). In particular, if $M$ is a closed manifold, $D_M$ being elliptic implies that $\text{dim}(ker D_M)$ is finite.

If $n=4k$, the spinor representation splits and there is a corresponding decomposision $S=S^+\oplus S^-$\label{def:spinor+-bundle}. The $Spin$ Dirac operator $D_M$ preserves this $\Z_2$-grading and exchanges the factors. We may restrict the $Spin$ Dirac operator to obtain operators $D^+_M:\Gamma(S^+)\rightarrow \Gamma(S^-)$ and $D^-_M:\Gamma(S^-)\rightarrow \Gamma(S^+)$ which satisfy $(D_M^+)^*=D_M^-$. We call $D_M^+$ the $Spin^+$ \emph{Dirac operator}\label{def:spin+diracoperator} of $M$.

If $W^{4k}$ is a compact spin manifold with boundary $M=\partial W$, then the restriction of $D_W^+$ to $M$ can be identified with the $Spin$ Dirac operator $D_M$ on $M$.

Let $\hat{A}$\label{def:Ahatgenus} denote the genus associated to the characteristic power series $(\sqrt{z}/2)/\sinh(\sqrt{z}/2)$ with corresponding multiplicative sequence $\{\hat{A}_k\}$ (see \cite[\textsection 19.]{MS74} and \cite[III \textsection 11]{LM89}). It is a power series in the Pontrjagin classes (or forms) of $M^{4k}$. In particular, for $k=1,2$ and $4$ we have

\begin{align}\label{eqn:Ahatgenus}
\hat{A}_1(p_1)&=-\frac{1}{24}p_1, \\
\hat{A}_2(p_1,p_2)&=\frac{1}{5760}\Big(-4p_2+7p^2_1\Big), \\
\hat{A}_4(p_1,p_2,p_3,p_4)&=\frac{1}{464486400}\Big(-192p_4+512p_3p_1+208p^2_2-904p_2p^2_1+381p^4_1\Big).
\end{align} If $x_i$ denote the formal roots of $TM$ for $i=1,...,2k$ (see \cite[p.9]{HBJ92} and \cite[III \textsection 11]{LM89}), then the Pontrjagin classes are given by the elementary symmetric functions\footnote{Recall for example that $\sigma_1(z_1,...,z_n)=\sum_{i=1}^n z_i$ and $\sigma_2(z_1,...,z_n)=\sum_{i<j}z_iz_j$.} in the square of the formal roots, i.e. $p_i(M)=\sigma_i(x^2_1,...,x^2_{2k})$, and the $\hat{A}$-genus is given by

\begin{equation}\label{eqn:Ahatgenusformalroots}
\hat{A}(M)=\prod_{i=1}^{2k}\frac{x_i/2}{\sinh(x_i/2)}.
\end{equation}

\begin{theorem}[$Spin$ Atiyah-Patodi-Singer index theorem]\label{thm:spinapsthm} \cite[Theorem (4.2)]{APSI75} Let $W^{4k}$ be a Riemannian spin manifold which is of product form near the boundary $M^{4k-1}=\partial W$. Assume that the $Spin^+$ Dirac operator $D^+_W$ and the restriction $D_M$ to the boundary $M$ satisfy the APS boundary condition (see \cite{APSI75}). Then the index of $D^+_W:\Gamma(S^+,P)\rightarrow \Gamma(S^-)$ is given by
$$\emph{index}(D^+_W)=\int_W \hat{A}(W)-\frac{h(D_M)+\eta(D_M)}{2}$$ where $\emph{index}(D^+_W)$, $h(D_M)$ and $\eta(D_M)$ are defined in Remark in Equations (\ref{def:index}),(\ref{def:kernel}) and (\ref{def:etainv}), $\hat{A}$ is the $\hat{A}$-genus in the Pontrjagin forms of the Riemannian metric on $W$.
\end{theorem}

If $E$ is a real oriented rank $2k$-vector bundle with formal splitting $E=E_1\oplus ... \oplus E_k$ into oriented 2-plane bundles and $y_j=e(E_j)$, we define

\begin{equation}\label{eqn:ahatpiofE}
\hat{A}_\pi(E):=\frac{1}{(2i)^k}\prod_{j=1}^k\frac{1}{\cosh(y_j/2)}.
\end{equation}

For a complex vector bundle $F$ with formal splitting $F=l_1\oplus ... \oplus l_k$ into complex line bundles bundles and $x_j=c_1(l_j)$, let

\begin{equation}\label{eqn:ahatthetaofF}
\hat{A}_\theta(F):=\frac{1}{2^k}\prod_{j=1}^k\frac{1}{\sinh(\frac{1}{2}(x_j+i\theta))}
\end{equation} where $0<\theta < \pi$.

Applying Donnelly's equivariant index theorem \cite[Theorem 1.2]{Do78} to the $Spin^+$ Dirac operator, we obtain the following.

\begin{theorem}[$Spin$ equivariant index theorem]\label{thm:equivindthmspin}
Let $W^{4k}$ be a compact Riemannian spin manifold which is of product form near the boundary $M^{4k-1}=\partial W$. Suppose that the $Spin^+$ Dirac operator $D^+_W$ and the restriction $D_M$ to the boundary $M$ satisfy the APS boundary condition. Let $g:W\rightarrow W$ be an isometry preserving the $Spin$ structure. Then the equivariant index of $D^+_W:\Gamma(S^+,P)\rightarrow \Gamma(S^-)$ is given by

$$
\emph{index}(D^+_W,g)=\sum_{N\subset W^g}a_{spin}(N)-\frac{\eta_g(D_M)+h_g(D_M)}{2},$$ where $\emph{index}(D^+_W,g)$, $h_g(D_M)$ and $\eta_g(D_M)$ are defined in Equations (\ref{def:equivindex}), (\ref{def:equivkernel}) and (\ref{eqn:equivetainv}), $N$ denotes a fixed point component of the action of $g$ on $W$ and $a_{spin}(N)$ is the corresponding local contribution.

If $N\subset W^g$ denotes a fixed point component without boundary, whose normal bundle splits as $\nu=\nu(\pi)\oplus \bigoplus_{0<\theta < \pi} \nu(\theta)$ in $W$ (see Remark \ref{rmk:normalbundledecomp}), then \cite[III. Theorem 14.11]{LM89}\label{def:localcontrspin}

\begin{equation}\label{eqn:localcontrspin}
a_{spin}(N)=(-1)^s\int_N\prod_{0<\theta \leq \pi}\hat{A}_\theta(\nu(\theta))\hat{A}(N),
\end{equation} where $s\in\{0,1\}$ depends on the action of $g$ on the $Spin$ structure (see \cite[III. Remark 14.12]{LM89}).
\end{theorem}

\subsection{Signature operator}\label{section:signatureoperator}

For a closed, oriented $4k$-dimensional manifold $M$, let $\text{sign}(M)$\label{def:signatureclosed} denote the signature of the non-degenerate quadratic form
\begin{equation}\label{eqn:quadformsign}
H^{2k}(M;\R)\times H^{2k}(M;\R)\rightarrow \R : (\alpha,\beta)\mapsto \int_M \alpha \cup \beta.
\end{equation} The integer $\text{sign}(M)$ is called the \emph{signature} of the closed manifold $M$.

Let $L$\label{def:Lgenus} be the genus associated to the characteristic power series $\sqrt{z}/\tanh(\sqrt{z})$ with corresponding multiplicative sequence $\{L_k\}$ (see \cite[\textsection 19.]{MS74} and \cite[III \textsection 11]{LM89}). It is a power series in the Pontrjagin classes (or forms) of $M$. In particular, for $k=1,2$ and $4$ we have

\begin{align}\label{eqn:Lgenus}
L_1(p_1)&=\frac{1}{3}p_1, \\
L_2(p_1,p_2)&=\frac{1}{45}\Big(7p_2-p^2_1\Big), \\
L_4(p_1,p_2,p_3,p_4)&=\frac{1}{14175}\Big(381p_4-71p_3p_1-19p^2_2+22p_2p^2_1-3p^4_1\Big).
\end{align}

Let $W^{2k}$ be an oriented, compact manifold with boundary $\partial M=W$. If $k$ is even, the \emph{signature} \label{def:signaturemfdwithboundary} of $W$ is defined as the signature of the quadratic form defined on $im(H^{k}(W,M;\R)\rightarrow H^k(W;\R))$ via the cup product. We will likewise denote it by $\text{sign}(W)$.

For the following definitions and discussion, see also \cite[p.63]{APSI75}.

Now suppose $W^{2k}$ is an oriented, compact Riemannian manifold which is of product form near the boundary $M^k=\partial W$. Let $\Omega^p(W)$ be the space of $p$-forms on $W$ and $\Omega(W):=\bigoplus_{p=0}^{2k}\Omega^p(W)$. Let $d:\Omega^p(W)\rightarrow \Omega^{p+1}(W)$ be the exterior derivative. The map $\tau:\Omega^*(W)\rightarrow \Omega^*(W)$ defined by $\tau(\omega)= i^{p(p-1)+k}\star \omega$ for $\omega \in \Omega^p(W)$ is an involution, where $\star$ is the Hodge star operator. Let $\Omega_\pm$ be the $\pm1$-eigenspaces of $\tau$ applied to $\Omega(W)$. Then $$A_W:=d+d^*=d-\star d\star:\Omega_+\rightarrow \Omega_-$$ is an elliptic operator which we will call the \emph{signature operator}. Since the metric is of product form near the boundary, we have $$A_W=\sigma\Big(\frac{\partial}{\partial u}+B_M\Big)$$ where $$B_M\omega=(-1)^{k+p+1}(\epsilon \star d-d\star)\omega$$ for $\epsilon =1$ if $\omega \in \Omega^{2p}(M)$ and $\epsilon=-1$ if $\omega \in \Omega^{2p-1}(M)$ (see \cite[p.63]{APSI75}). The operator $B_M$ is formally self-adjoint and preserves the parity of the forms, so that there is a decomposition $B_M=B_M^{ev}\oplus B_M^{odd}$\label{def:Bevsignatureoperator}. The operator $B^{ev}_M$ is sometimes called the \emph{odd signature operator} of $M$.

\begin{theorem}[Signature Atiyah-Patodi-Singer index theorem]\label{thm:signAPSindthm} \cite[Theorem (4$\cdot$14)]{APSI75}
Let $W^{4l}$ be a compact oriented Riemannian manifold which is of product form near the boundary $M=\partial W$. Then
$$\emph{sign}(W)=\int_W L(W)-\eta(B^{ev}_M)$$ where $L(W)$ is the $L$-genus in the Pontrjagin forms of $W$ and $\eta(B^{ev}_M)$ by Equation (\ref{def:etainv}).
\end{theorem}

For the following definitions and discussion, see also \cite[pp.408-409]{APSII75}.

Now let $W^{2k}$ be a compact, oriented, Riemannian manifold which is of product form near the boundary $M^{2k-1}=\partial W$, $A_W$ the signature operator and $B_M=B^{ev}_M\oplus B^{odd}_M$ the operator from above. 

Assume that $G$ is a compact Lie group acting via orientation preserving isometries on $W^{2k}$. The induced action of $G$ on sections commutes with the operator $A_W$. Furthermore, $G$ acts on $\tilde{H}^k(W;\R)$ and preserves the quadratic form defining the signature of a manifold with boundary (see above) which is symmetric if $k$ is even and skew-symmetric if $k$ is odd. Complexify and consider the corresponding hermitian form. Now, any $G$-invariant inner product on $\tilde{H}^k(W;\R)$ will induce a $G$-invariant decomposition $\tilde{H}^k(W;\R)=\tilde{H}^k_+\oplus \tilde{H}^k_-$ such that the hermitian form is positive definite on $\tilde{H}^k_+$ and negative definite on $\tilde{H}^k_-$. The virtual representation $\text{sign}(G,W):=\tilde{H}^k_+-\tilde{H}^k_-$ is called the \emph{$G$-signature} of $W$ and \label{def:equivariantsignature}
\begin{equation}\label{eqn:equivsign}
\text{sign}(g,W):=tr(g|_{\tilde{H}^k_+})-tr(g|_{\tilde{H}^k_-})
\end{equation} the \emph{equivariant signature} of $W$ with respect to $g\in G$.

\begin{theorem}[Signature equivariant index theorem]\label{thm:equivindexthmsign}
\cite[Theorem 2.1]{Do78} Let $W^{2k}$, $M^{2k-1}$, $G$, $A_W$ and $B_M=B^{ev}_M\oplus B^{odd}_M$ be as above and suppose they satisfy the APS boundary condition. Then for each $g\in G$,

\begin{equation}
\emph{sign}(g,W)=\sum_{N\subset W^g}a_{sign}(N)-\eta_g(B^{ev}_M)
\end{equation} where $\eta_g(B^{ev}_M)$ is defined using Equation (\ref{eqn:equivetainv}). 

If $N\subset W^g$ denotes a fixed point component without boundary, whose normal bundle splits as $\nu=\nu(\pi)\oplus \bigoplus_{0<\theta < \pi} \nu(\theta)$ in $W$ (see Remark \ref{rmk:normalbundledecomp}), then \cite[III. Theorem 14.5]{LM89}\label{def:localcontributionsignature}

\begin{equation}\label{eqn:localcontrsign}
a_{sign}(N)= \int_N\prod_{0<\theta\leq \pi}L_\theta(\nu(\theta))L(N),
\end{equation} where for any oriented real vector bundle $E$,
\begin{equation}
L_\pi(E):=e(E)(L(E))^{-1},
\end{equation} and for any complex vector bundle $F$ with formal splitting $F=l_1\oplus ... \oplus l_k$ into complex line bundles with $x_j=c_1(l_j)$,
\begin{equation}
L_\theta(F):=\prod_{j=1}^k\coth\big(x_j+\frac{i\theta}{2}\big),
\end{equation} for $0<\theta <\pi$.
\end{theorem}

\subsection{Eta-invariant of a covering}

Let $M^{2n+1}$ be a closed, oriented, Riemannian manifold and let $\pi:\tilde{M}\rightarrow M$ be a regular covering with finite covering group $G$. The metric on $M$ lifts to a metric on $\tilde{M}$ and any elliptic self-adjoint operator $D_M:\Gamma(E)\rightarrow \Gamma(F)$ (where $E$ and $F$ are vector bundles on $M$ with a smooth inner product) lifts to an elliptic self-adjoint operator $D_{\tilde{M}}:\Gamma(\pi^*E)\rightarrow \Gamma(\pi^*F)$ which is equivariant with respect to the action of $G$ by deck transformation. For each irreducible unitary representation $\alpha:G\rightarrow U(k)$, there is a flat vector bundle $E_\alpha:=\tilde{M}\times_\alpha \C^k\rightarrow M$ and a twisted operator $D_{M,E_\alpha}:\Gamma(E)\otimes E_\alpha\rightarrow \Gamma(F)\otimes E_\alpha$.

\begin{theorem}\label{thm:etainvcov}
Let $M$, $\tilde{M}$, $G$, $D_{\tilde{M}}$ and $D_{M,E_\alpha}$ be as above. Then 
\begin{equation}
\eta(D_{M,E_\alpha})=\frac{1}{|G|}\sum_{g\in G}\eta_g(D_{\tilde{M}})\cdot\chi_\alpha(g)
\end{equation} where $\chi_\alpha$ is the character of $\alpha$, $\eta(D_{M,E_\alpha})$ is defined via Equation (\ref{def:etainv}) and $\eta_g(D_{\tilde{M}})$ via Equation (\ref{eqn:equivetainv}).
\begin{proof}
The statement in general follows from \cite[(2$\cdot$14)]{APSII75} and the orthogonality relations of irreducible characters. For the special case of the signature operator, see also \cite[Theorem 3.4.]{Do78}.
\end{proof}

\end{theorem}

\subsection{Applications to positive scalar curvature}\label{subsec:posscalreletainv}

\subsubsection{Vanishing theorems}

For the following result, which is due to Lichnerowicz \cite{Li63}, see also \cite[II. Corollary 8.9]{LM89} and \cite[II. Theorem 8.11]{LM89}.

\begin{theorem}[$Spin$ vanishing theorem]\label{thm:spinvanishingtheorem}
Let $M^{4k}$ be a closed spin manifold and $D^+_M$ its $Spin^+$ Dirac operator. If the Riemannian metric on $M$ has non-negative scalar curvature everywhere and positive scalar curvature at some point, then $ker D_M=0$ and consequently $\text{index}(D^+_M)=\hat{A}(M)=0$.
\end{theorem}

There also is a version for manifolds with boundary (see \cite[Theorem (3$\cdot$9)]{APSII75}).

\begin{theorem}[$Spin$ vanishing theorem with boundary]\label{thm:spinvanishingthmboundary}
Let $W^{2l}$ be a compact spin manifold with boundary $M^{2l-1}$. Let $D^+_W$ be the $Spin^+$ Dirac operator on $W$ and $D_M$ the $Spin$ Dirac operator on the boundary. If there is a Riemannian metric on $W$ which is of product form near the boundary and which has non-negative scalar curvature everywhere and positive scalar curvature at some point on $M$, then $\emph{index}(D^+_W)=0$ and $\emph{ker}(D_M)=0$.

\end{theorem}

\section{Appendix B. C++ code}

Unfortunately, it is unknown to the author whether there exist number theory methods to solve the arithmetics of Proposition \ref{prop:eellskuipershimadaquodiffvalues}. As a resort, the following C++ code counts the number of different values of the Eells-Kuiper invariant of the Shimada projective spaces.

\begin{verbatim}
#include <iostream>
#include <iomanip> // for setw
using namespace std;

int main() {

int counter;
counter=0;
int countermu, countermuquo, helpcountermu, helpcountermuquo;
countermu=0;
countermuquo=0;
int n, nn, m;
n=16255;

for (int i=0;i<n;i++){
    helpcountermu=0;
    helpcountermuquo=0;
 for (int k = i; k < n; k++) {

     int mui, muk, a, b, c, d;
     int muiquoplus, mukquoplus, muiquominus, mukquominus;

    mui=i*(i+1)%16256;
    muk=k*(k+1)%16256;

    a=2*i*(i+1)+127*(2*i+1);
    muiquoplus=a%65024;
    b=65024+2*i*(i+1)-127*(2*i+1);
    muiquominus=b%65024;

    c=2*k*(k+1)+127*(2*k+1);
    mukquoplus=c%65024;
    d=65024+2*k*(k+1)-127*(2*k+1);
    mukquominus=d%65024;

if (mui==muk) {
  if (k!=i) {
  
//If the Eells-Kuiper invariants of the Shimada sphere are equal,
//we don't count it as a "new" distinct value
//Thus, we increment helpcountermu and if it is non-zero,
//we don't increment countermu

      helpcountermu++;

//If the Eells-Kuiper invariants of the quotients are equal,
//we don't count it as a "new" distinct value
//Thus, we increment helpcountermu and if it is non-zero,
//we don't increment countermuquo

if (muiquoplus==mukquoplus && muiquominus==mukquominus)
  {helpcountermuquo++;}
else if (muiquoplus==mukquominus && muiquominus==mukquoplus)
  {helpcountermuquo++;}
    }
  }
}
 if (helpcountermu==0) {countermu++; }
 if (helpcountermuquo==0) {countermuquo++;}
}


cout << "Number mu values (spheres): " << countermu << endl;
cout << "Number of mu values (quotients): "  << countermuquo << endl;
return 0;
}

\end{verbatim}


\printbibliography

\end{document}